\documentclass{article}

\newtheorem{theorem}{Theorem}[section]
\newtheorem{definition}{Definition}

\begin{document}

\title{On a stellar structure for a stellar manifold}

\author{Sergey Nikitin \\
Department of Mathematics\\
Arizona State University\\
Tempe, AZ 85287-1804\\
 nikitin@asu.edu \\
 http://lagrange.la.asu.edu}

\maketitle

{\bf Abstract} It is well known that a compact two dimensional surface is homeomorphic to a polygon with the edges identified in pairs. This paper not only presents a new proof of this statement but also generalizes it for any connected $n$-dimensional stellar manifold with a finite number of vertices.

\vspace{0.1cm}

\section{Introduction}

\vspace{0.1cm}

It is well known (see e.g. \cite{Massey}) that a compact $2$-dimensional surface is homeomorphic to a polygon with the edges identified in pairs.
 In this paper we show that any connected $n$-dimensional stellar manifold $M$ with a finite number of vertices possesses a similar property: $M$ is stellar equivalent to 
$$
       a\star (S/\simeq),
$$
where $a \notin S$ is a vertex, $S$ is $(n-1)$-dimensional stellar sphere. The equivalence relation "$\simeq$" satisfies the following conditions:
\begin{itemize}
\item[] No $(n-1)$-dimensional simplex of $S$ has two vertices that are equivalent to each other.
\item[] For any $(n-1)$-dimensional simplex $g$ of $S$ there might exist not more than one (different from $g$) $(n-1)$-dimensional simplex $p$ in $S$ such that any vertex of $g$ is equivalent to some vertex of $p.$
\end{itemize}
$a \star (S/\simeq)$ is not any more a stellar manifold in its usual sense.
We take the ball $a\star S$ and identify equivalent generators from $S$ but,
at the same time, we distinguish any two generators $a\star p$ and $a\star g$ even when $p\simeq g.$ In the sequel, 
it is convenient to use the notation $a\star (S/\simeq)$ for the resulted manifold even if it is slightly misleading.

\section{Main results}
We begin with recalling the basic definitions of stellar theory \cite{Glaser}, \cite{Lickorish}.
A stellar $n$-manifold $M$ can be identified with the sum of its $n$-dimensional simplexes ($n$-simplexes):
$$
M=\sum_{i=1}^n g_i
$$
with coefficients from ${\rm Z}_2.$ We will call $\{g_i\}_{i=1}^n$ generators of $M.$

All vertices in $M$ can be enumerated and any $n$-simplex $s$ from $M$ corresponds to the set
of its vertices
$$
s=(i_1 \; i_2 \; \dots \; i_{n+1}),
$$
where $i_1 \; i_2 \; \dots \; i_{n+1}$ are integers. 

The boundary operator $\partial $ is defined on a simplex as
$$
\partial (i_1 \; i_2 \; \dots \; i_{n+1}) = (i_1 \; i_2 \; \dots \; i_n ) + (i_1 \; i_2 \; \dots \; i_{n-1} \; i_{n+1} ) + \dots + (i_2 \; i_3 \; \dots \; i_{n+1})
$$
and linearly extended to any complex, i.e.
$$
\partial M = \sum_{i=1}^n \partial g_i.
$$
A manifold is called closed if $\partial M =0.$

If two simplexes $(i_1 \; i_2 \; \dots i_m) $ and $(j_1 \; j_2 \; \dots j_n)$ do not have common vertices then one can define their join
$$
(i_1 \; i_2 \; \dots i_m) \star (j_1 \; j_2 \; \dots j_n)
$$
as the union
$$
(i_1 \; i_2 \; \dots i_m) \cup (j_1 \; j_2 \; \dots j_n).
$$

\vspace{0.1cm}

If two complexes $K=\sum_i q_i $ and $ L = \sum_j p_j$ do not have common vertices then 
their join is defined as
$$
K\star L = \sum_{i,j} q_i \star p_j.
$$

If $A$ is a simplex in a complex $K$ then we can introduce its link:
$$
lk(A,K) = \{ B \in K \; ; \; A \star B \in K \}.
$$
The star of $A$ in $K$ is $A \star lk(A,K).$ Thus,
$$
K = A\star lk(A,K) + Q(A,K),
$$
where the complex $Q(A,K)$ is composed of all the generators of $K$ that do not contain $A.$ A complex with generators of the same dimension is called a uniform complex.

\begin{definition} ({\bf Subdivision})
Let $A$ be a simplex of a complex $K.$ Then any integer $a$ which is not a vertex of $K$ defines starring of
$$
K=A\star lk(A,K) + Q(A,K)
$$
at $a$ as
$$
 \hat{K}=a\star \partial A \star lk(A,K) + Q(A,K).
$$
This is denoted as
$$
\hat{K} =  (A \; a)K.
$$
\end{definition}

The next operation is the inverse of subdivision. It is called a stellar weld and defined as follows.

\begin{definition} ({\bf Weld})
  Consider a complex
  $$
        \hat{K}=a\star  lk(a,\hat{K}) + Q(a,\hat{K}),
  $$
with $lk(a,\hat{K}) = \partial A \star B$ where
$B$ is a subcomplex in $\hat{K},$  $A$ is a simplex and $A\notin \hat{K} .$
Then the (stellar) weld $(A\; a)^{-1} \hat{K}$ is defined as
  $$
       (A\; a)^{-1} \hat{K} =  A \star B  + Q(a,\hat{K}).
  $$
\end{definition}

A stellar move is one of the following operations: subdivision, weld, enumeration change on the set of vertices.
Two complexes $M$ and $L$  are called stellar equivalent if one is obtained from the other by a finite sequence of stellar moves.
It is denoted as $M \sim L.$ We also say that $M$ admits triangulation $L.$\\

If a complex $L$ is stellar equivalent to $(1 \; 2 \; \dots \; n+1)$ then $L$ is called a stellar 
$n$-ball. On the other hand, if $K \sim \partial (1 \; 2 \; \dots \; n+2)$ then $K$ is a stellar
$n$-sphere.

\begin{definition} ({\bf Stellar manifold})
\label{stellar_def}
Let $M$ be a complex. If, for every vertex $i$ of $M,$ the link $lk(i,M)$ is either a stellar $(n-1)$-ball or a stellar $(n-1)$-sphere, then $M$ is  a
stellar $n$-dimensional manifold ($n$-manifold).
\end{definition}

If $i$ is a vertex of $M$ then
$$
M= i\star lk(i,M) + Q(i,M).
$$
If $\partial M=0,$ then $Q(i,M)$ is a stellar manifold.\\
Indeed, consider an arbitrary vertex $j$ of $Q(i,M).$ Then
$$
lk(i,M)=j\star lk(j,lk(i,M)) + Q(j,lk(i,M))
$$
and
$$
Q(i,M)=j\star lk(j,Q(i,M)) + Q(j,Q(i,M)).
$$
Since $M$ is a stellar manifold and $\partial M = 0$
$$
i\star lk(j,lk(i,M)) + lk(j,Q(i,M)) 
$$
is a stellar sphere. Hence, it follows from \cite{Newman} that $lk(j,Q(i,M))$ is either a stellar ball or a stellar sphere.\\
In the sequel it is convenient to consider an equivalence relation on the set of vertices of a stellar manifold. Among all possible such 
equivalence relations we are mostly interested in those that meet certain regularity properties underlined by the following definition.
\begin{definition} ({\bf Regular equivalence})
Given a stellar manifold $M,$ an equivalence relation "$\simeq$" on the set of vertices from $M$ is called regular if it meets the following conditions:
\begin{itemize}
\item[(i)] No generator $g \in M$ has two vertices that are equivalent to each other.
\item[(ii)]For any generator $g\in M$ there might exist not more than one generator $p\in M \setminus g$ such that any vertex of $g$ is either equal or equivalent to some vertex of $p.$ We call such two generators equivalent, $g \simeq p.$ 
\end{itemize}
\end{definition}

Throughout the paper we use notation $a\star (S/\simeq)$ to identify any two equivalent generators from $S$ but distinguish the corresponding generators from $a\star S.$
The next result of the paper is  valid for any connected stellar manifold with a finite number of generators.

\begin{theorem}
\label{triangulation}
A connected stellar $n$-manifold $M$ with a finite number of generators admits a triangulation
$$
N = a\star (S/\simeq) ,
$$
where $a\notin S$ is a vertex, $S$ is a stellar $(n-1)$-sphere and "$\simeq $" is a regular equivalence relation. Moreover, if $M$ is closed then for any generator $g\in S$ there exists exactly one generator $p\in S\setminus g$ such that $g \simeq p.$ 
\end{theorem}
{\bf Proof.}
Let us choose an arbitrary generator $g \in M$ and an integer $a$ that is not a vertex of $M.$ Then
$$
M \sim (g \; a ) M  \mbox{  and  }
(g \; a ) M =a \star \partial g + M\setminus g,
$$
where $M\setminus g$ is defined by all the generators of $M$ excluding $g.$
We construct $N$ in a finite number of steps.
Let $N_0=(g \; a ) M.$ Suppose we constructed already $N_k$ and there exists a generator $p\in  Q(a,N_k)$ that has at least one common $(n-1)$-simplex with $lk(a, N_k).$
Without loss of generality, we can assume that 
$$
 p = (1 \; 2\; \dots \; n+1).
$$
and $(1 \; 2\; \dots \; n)$ belongs to $lk(a, N_k).$
If the vertex $(n+1)$ does not belong to  $lk(a, N_k)$ then
$$
N_{k+1} = ((a\; n+1 )\; b)^{-1}((1 \; 2 \; \dots \; n) \; b) N_k,
$$
where $b\notin N_k.$
If the vertex $(n+1)$ belongs to $lk(a, N_k)$ then after introducing a new vertex $d\notin N_k$ we take
$$
L = ((a\; d )\; b)^{-1}((1 \; 2 \; \dots \; n) \; b)(N_k\setminus p + (1\;2\;\dots \; n \;d)),
$$
where $b\notin (N_k\setminus p + (1\;2\;\dots \; n \;d )),$ and
$$
N_{k+1} =  a \star (lk(a,L)/\simeq ) + Q(a,L)
$$
endowed with the equivalence $d \simeq (n+1).$\\
By construction
$$
N_{k+1} = a \star (lk(a,N_k) +\partial p ) + Q(a,N_k) \setminus p
$$
if $(n+1)\notin lk(a, N_k).$ Otherwise,
$$
N_{k+1} = a \star ((lk(a,N_k) +\partial g)/\simeq ) + Q(a,N_k) \setminus p,
$$
where $g=(1 \; 2\; \dots \; n \; d),\;\;d\simeq (n+1).$

Since $M$ is connected and has a finite number of generators there exists a natural number $m$ such that
$$
N_m =  a \star (S/\simeq)
$$
where $S$ is a stellar $(n-1)$-sphere and "$\simeq $" is a regular equivalence relation.\\
If $M$ is closed, then $\partial N_m = 0,$ and therefore, for any generator $g \in S$ there exists exactly one generator $p\in S \setminus g$ such that $g \simeq p.$\\
Q.E.D.\\

When working with stellar manifolds we use topological terminology in the following sense. We say that a stellar $n$-manifold $M$ has a certain topological property if its standard realization in ${\rm R}^{2n+1}$ (see e.g \cite{Alexandrov}) has this property.\\
Let $\pi (M)$ denote fundamental group of a manifold $M.$
Then Theorem \ref{triangulation} together with Seifert -- Van Kampen theorem (see e.g. \cite{Massey}) lead us to the next result.

\begin{theorem}
\label{fundamental}
If $M$ is a connected stellar $n$-manifold and $n>2$ then
$$
\pi (M) = \pi (S/\simeq ),
$$
where $S$ is a stellar $(n-1)$-sphere and "$\simeq $" is a regular equivalence relation.
\end{theorem}
{\bf Proof.}
By Theorem \ref{triangulation}  $M$ admits a triangulation
$$
N = a\star (S/\simeq ),
$$
where $a\notin S$ is a vertex, $S$ is a stellar $(n-1)$-sphere and "$\simeq $" is a regular equivalence relation.\\
Using the triangulation $N$ we can define open subsets $U$ and $V$ in $M$ as follows
\begin{eqnarray*}
U &=& a\star (S/\simeq )\setminus a \\
V &=& a\star (S/\simeq ) \setminus (S/\simeq ).
\end{eqnarray*}
Clearly,
$$
M = U \cup V
$$
and $V$ is homeomorphic to an interior of the $n$-ball. Thus, $V$ is simply connected and by Seifert -- Van Kampen theorem there exists an epimorphism
$$
\psi \;:\; \pi(U)\; \longrightarrow \; \pi(M)
$$
induced by inclusion $U \subset M .$ The kernel of $\psi $ is the smallest normal subgroup containing the image of the homomorphism
$$
\varphi\;:\;\pi(U\cap V) \; \longrightarrow \; \pi(U)
$$
induced by inclusion $U\cap V \subset U.$\\

The open set $U\cap V $ is homeomorphic to 
$$
0 < x_1^2+ x_2^2 +\dots +x_n^2 <1,\;\;\mbox{ where } \;\; (x_1,\; x_2 ,\; \dots ,\;x_n) \in {\rm R}^n.
$$
The fundamental group of this manifold is trivial for $n>2.$ Hence, the kernel of the epimorphism $\psi $ is trivial and
$$
\pi(M) = \pi(U).
$$
On the other hand, $S/\simeq $ is a deformation retract of $U,$ and therefore,
$$
 \pi(U) = \pi (S/\simeq).
$$
Q.E.D.\\

\bibliographystyle{amsplain}

\end{document}